\documentclass[5p,times]{elsarticle}

\usepackage{amssymb}
\usepackage{amsthm}
\usepackage{makecell}
\usepackage{amsmath}
\usepackage{color,xcolor}
\newtheorem{remark}{Remark}
\usepackage[utf8]{inputenc}
\usepackage{booktabs, caption, makecell}

\usepackage{threeparttable}
\usepackage{epigraph}
\usepackage{enumitem}
\usepackage{hyperref}
\usepackage{subfigure}
\usepackage{mathtools}
\usepackage{url}
\usepackage{color}
\usepackage{bbm}
\usepackage{amsmath, amssymb, amsfonts, amsthm}
\usepackage{siunitx}
\usepackage[linesnumbered, lined, boxed]{algorithm2e}

\usepackage{multirow}
\usepackage[figuresright]{rotating}
\usepackage{amsmath}

\usepackage{tablefootnote}

\begin{document}
\begin{frontmatter}

\title{Symbolic Discovery of Iterative Algorithms: \\ A Continuous Latent Space Bayesian Optimization Framework}
\author{Tongjia Liu}
\author{Ilias Mitrai\fnref{label2}}
\fntext[label2]{Corresponding author}
\ead{imitrai@che.utexas.edu}

\address{McKetta Department of Chemical Engineering, The University of Texas at Austin, Austin, TX 78712, US}

\begin{abstract}
In this paper, we consider the automated discovery of iterative optimization algorithms. We formulate the algorithm discovery task as a discrete optimization problem and search for new update functions using latent space Bayesian Optimization. The proposed framework first learns a continuous representation of the discrete space of update functions using variational autoencoders, transforming the algorithm discovery task from a discrete to a continuous search problem. The continuous representation is subsequently used to search for new algorithms using Bayesian optimization. Application to two case studies shows that the proposed approach can discover new update functions in symbolic form without any assumptions on the functional form of the update function. Moreover, the computational time required to discover the new update functions is lower than existing mathematical programming-based approaches.  
\end{abstract}


\end{frontmatter}

\section{Introduction}
Iterative algorithms are widely used to solve computational tasks arising in chemical engineering. Given a computational task, denoted as $F(x)=0$, iterative algorithms use an update function $g$ to update the estimate of the solution, $x_{i+1}=g(x_i)$, until a predetermined convergence criterion is met. Update functions are usually developed using domain knowledge, which is inherently time-consuming. This limitation has motivated the development of automated algorithm discovery approaches. 

Existing approaches differ in the construction, i.e., parameterization of update functions, and navigation of the search space. For example, neural networks have been used to parameterize the update function \cite{andrychowicz2016learning, chen2022learning}. Despite promising results, this approach results in black-box update functions whose convergence properties are difficult to analyze. Recently, basis functions have been used to represent $g$ \cite{mitsos2018optimal}. While this approach can yield interpretable update functions, the functional form is predetermined, limiting the space of algorithms that can be explored. 

An alternative is to represent the update function as an expression tree and search using evolutionary algorithms, such as genetic programming \cite{chen2024symbolic} and large language models \cite{romera2024mathematical}. These methods yield interpretable update functions without any a priori assumptions regarding the functional form of $g$. However, the search is stochastic and computationally inefficient. 
\begin{figure*}
    \centering
    \includegraphics[trim = 0 0 0 0, clip, scale=0.53]{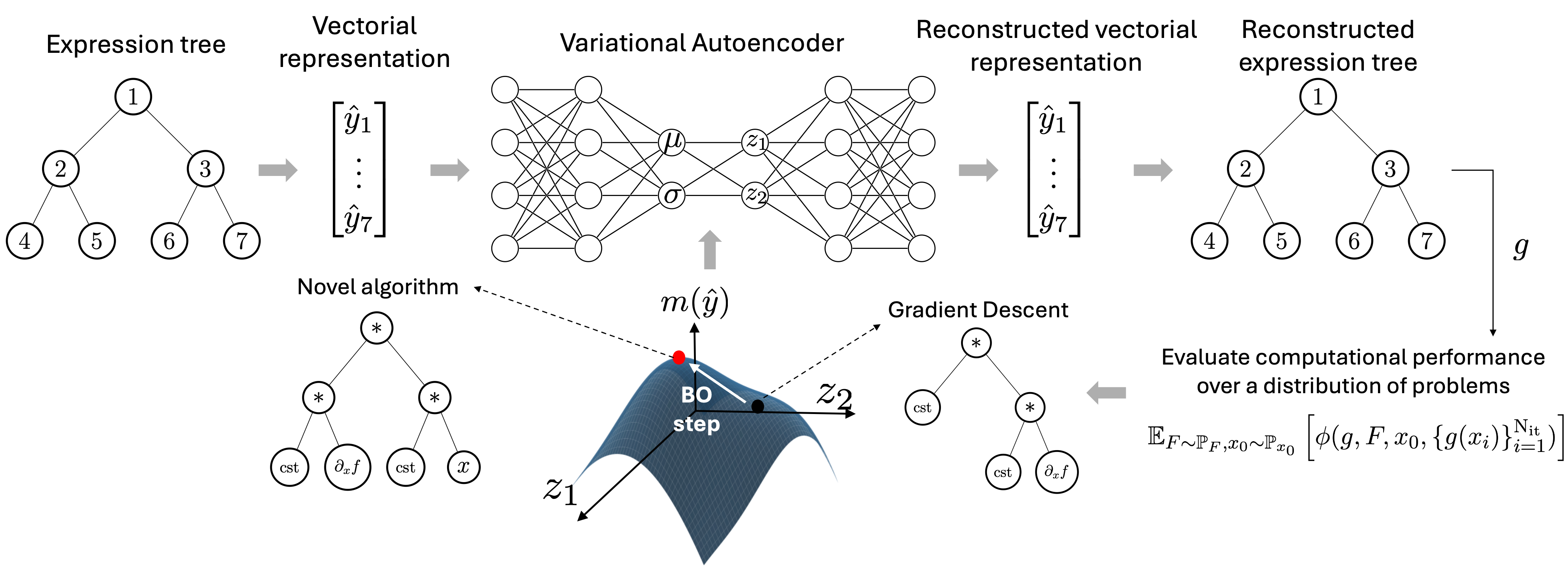}
    \caption{Proposed Variational Autoencoder-based Bayesian Optimization search for iterative algorithms}
    \label{fig:data driven ades framework}
\end{figure*}
In recent work, we have proposed a mixed-integer optimization approach to search the space of update functions represented by expression trees \cite{mitrai2025automated}. This approach combines the interpretability of expression trees with performance guarantees from global optimization. However, the search is computationally expensive and cannot be used if the algorithm evaluation metric does not have a closed-form expression, e.g., solution time. Overall, existing automated discovery methods face a trade-off between interpretability and scalability.

In this work, we propose a Bayesian Optimization approach to discover update functions in symbolic form. The proposed framework is presented in Figure~\ref{fig:data driven ades framework} and is based on the observation that if the update function is fixed, the algorithm discovery task reduces to an evaluation task, i.e., evaluating the algorithm's computational performance. However, this evaluation can be computationally expensive, especially for high-dimensional computational tasks. This structure lends itself to the application of BO, a sample-efficient search strategy well-suited to problems with expensive function evaluations. 

The proposed framework has two components. First, we learn a continuous representation of the discrete space of update functions, represented as expression trees, using Variational Autoencoders (VAEs) \cite{kingma2013auto}. This learning step transforms the discrete algorithm discovery task into a continuous one. This latent space is used to search for new update functions using Bayesian Optimization. We use the proposed approach to search for new update functions for solving unconstrained optimization problems. The results show that the proposed approach can successfully discover new update functions with lower computational effort than existing mathematical programming approaches.

\section{The Algorithm Discovery Problem}
We consider the case where an update function $g$ is used to solve unconstrained optimization problems $\min_{\textbf{x}} F(\textbf{x})$, $\textbf{x} \in \mathbb{R}^{\rm{N_{x}}}$, over a distribution of objective functions $F$ and starting points $\textbf{x}_{0}$, denoted as $\mathbb{P}_{F}$ and $\mathbb{P}_{x_{0}}$, respectively. We assume that all the problems have the same number of variables, $\rm{N_{x}}$, with $\mathcal{N}_{x}$ being the set of variables. We define as $\mathcal{I} = \{1,\dots,\mathrm{N_{it}}\}$ the set of iterations and $\mathrm{N_{it}}$ as the maximum number of iterations. 

We want to discover an update function $g$ that minimizes the expected value of a performance measure $\phi$, e.g., the norm of the gradient at the last iteration, over the distribution of $F$ and $\textbf{x}_{0}$ for $\mathrm{N_{it}}$ iterations. This task corresponds to the following algorithm discovery problem
\begin{equation}
    \begin{aligned}
\min_{g} \ & \mathbb{E}_{F\sim \mathbb{P}_{F}, \textbf{x}_0\sim \mathbb{P}_{x_0}} \left[ \phi \big(g, F, \textbf{x}_0, \{\textbf{x}_{i}(g,F,\textbf{x}_{0})\}_{i=1}^{\mathrm{N_{it}}} \big) \right],
\end{aligned}
\end{equation}
where $\{\textbf{x}_{i}(g, F,x_{0})\}_{i=1}^{\mathrm{N_{it}}}$ represents the trajectory of iterates generated by update function $g$ for a given objective $F$ and starting point $\textbf{x}_{0}$. The iterative scheme for variable $k$ is given by the following equation
\begin{equation}
    x_{k,i} = x_{k,i-1} + g(x_{k,i-1}, \partial_{x_{k}} F(\mathbf{x}_{i-1})) \ \forall k \in \mathcal{N}_x, \forall i \in \mathcal{I}.
\end{equation}
This problem can be transformed into a deterministic one by sampling $\mathrm{N}_{F}$ functions and $\mathrm{N_{x_{0}}}$ starting points, leading to the following deterministic algorithm discovery task
\begin{equation}
    \begin{aligned}
\min_{g} \ & \frac{1}{\mathrm{N_{F}} \ \mathrm{N_{x_{0}}}} \sum_{f=1}^{\mathrm{N_{F}}} \sum_{p=1}^{\mathrm{N_{x_{0}}}} \phi \big(g, F_{f}, \textbf{x}_{0}^{fp}, \{\textbf{x}_{fpi}\}_{i=1}^{\mathrm{N_{it}}}   \big),
\end{aligned}
\end{equation}
where $\textbf{x}_{fpi}$ is the vector of variables for function $f$, starting point $p$ and iteration $i$ and the associated entry for variable $k$ is denoted as $x_{fpki}$. $\textbf{x}_{0}^{fp}$ represents the values of the variables for starting point $p$ and function $f$. The formulation presented above is generic and can be readily adapted to discover iterative algorithms for any computational task, e.g., solving systems of algebraic equations or partial differential equations. The solution of the algorithm discovery task presents two challenges. The first concerns constructing the space of all possible update functions, and the second concerns the search technique used to navigate the search space.

\subsection{Symbolic representation of update functions}
We represent $g$ as a binary expression tree of depth $d$ with $N = 2^{d+1}-1$ nodes. We define $\hat{\mathcal{N}} = \{1, \dots, N\}$ as the set of nodes and $\mathcal{T}=\{2^{d},\dots,2^{d+1}-1\}$ as the set of terminal nodes, i.e., leaf nodes. Each internal node $n \in \hat{\mathcal{N}}\setminus \mathcal{T}=\mathcal{N}_{int}$ has two children, nodes $2n$ and $2n+1$. At each node, a binary $\mathcal{B} = \{+, -, *,/\}$ and unary operator $\mathcal{U} = \{\sin, \cos, \exp\}$ can be assigned as well as an operand $\mathcal{L} = \{x_k, \partial_{x_k} F, \text{cst}\}$, where $\text{cst}$ represents the constants of the update function. We denote $\mathcal{O} = \mathcal{B} \cup \mathcal{U} \cup \mathcal{L}$ as the set of all possible operators and operands that can be assigned to the nodes.

We define a binary variable $y_{on}$ which is equal to one if operator $o$ is assigned to node $n$ and zero otherwise. We also define $c_{n}$ as the value of the constant assigned to node $n$. These variables provide an exact parameterization of the update function, since they can explicitly determine the functional form and associated parameters. This expression tree representation can represent any functional form, provided that the tree depth and the operations considered are sufficient \cite{smith1999symbolic}. 

The assignment of operators to the nodes of the expression tree must follow arity rules. These constraints are adapted from \cite{cozad2018global,kim2023learning, mitrai2026constrained}. The first set of constraints guarantees that at most one operator can be assigned per node as follows
\begin{equation}\label{eq: one operator}
    \sum_{o \in \mathcal{O}} y_{on} \leq 1 \ \ \forall n \in \hat{\mathcal{N}}.
\end{equation}
The assignment of an operator at node $n$ depends on the operator of the parent node and affects the operators assigned to descendant nodes. For example, a binary operator can be assigned only to internal nodes. These constraints are enforced via the following equations
\begin{equation} \label{eq: operator assignment 1}
    \begin{aligned}
        \sum_{o \in \mathcal{B} \cup \mathcal{U}} y_{on} & = \sum_{o \in \mathcal{O}} y_{o,2n+1} \ \ \forall n \in \mathcal{N}_{\rm{int}}
    \end{aligned}
\end{equation}
\begin{equation} \label{eq: operator assignment 2}
    \begin{aligned}
        \sum_{o \in \mathcal{B}} y_{on} & = \sum_{o \in \mathcal{O}} y_{o,2n} \ \ \forall n \in n \in \mathcal{N}_{\rm{int}}.
    \end{aligned}
\end{equation}
Eq.~\ref{eq: operator assignment 1} enforces that if a binary or unary operator is assigned to a node $n$, then one operator is assigned to node $2n+1$. Eq.~\ref{eq: operator assignment 2} enforces that if a binary operator is assigned to node $n$ then an operator is assigned to node $2n$. 

Given the expression tree representation of the update function, we obtain the following algorithm discovery task, which is a discrete optimization problem
\begin{equation}
    \begin{aligned}
\min_{y,c} \ & \frac{1}{\mathrm{N_{F}} \ \mathrm{N_{x_{0}}}} \sum_{f=1}^{\mathrm{N_{F}}} \sum_{p=1}^{\mathrm{N_{x_{0}}}} \phi \big(y,c,F_{f}, \textbf{x}_{0}^{fp}, \{\textbf{x}_{fpi}\}_{i=1}^{\mathrm{N_{it}}}   \big)\\
\text{s.t.} \ \ & \text{Eq.}~4,5,6\\
& x_{fpk,i} = x_{fpk,i-1} + g(y,c,x_{fpk,i-1}, \partial_{x_{k}} F_{f}(\mathbf{x}_{fpi-1})) \\
& \ \forall k \in \mathcal{N}_x, \forall i \in \mathcal{I}, f \in \mathcal{F}, p \in \mathcal{P},
\end{aligned}
\end{equation}
where $y$, $c$ represent all the binary and continuous variables, respectively. $\mathcal{F}$ is the set of objective functions and $\mathcal{P}$ the set of starting points considered.

\section{Data-Driven Algorithm Discovery}
We observe that if the binary variables $y$ are fixed, the resulting problem is a nonlinear programming problem that seeks the values of the constants, $c$, that optimize the expected computational performance. Motivated by this structure, we propose a data-driven solution approach. Specifically, we use BO to search for the functional form of the update function and treat the search for the constants that minimize the expected computational performance as a black box.

\subsection{Continuous representation of iterative algorithms via Variational Autoencoders}
In this section, we present the transformation of the discrete search space into a continuous one using VAEs (see top part of Figure~\ref{fig:data driven ades framework}). First, the binary-variable representation of the expression tree is converted into an integer representation. We order the set of operators $o \in \mathcal{O}$ and convert the binary variables for each node $n \in \hat{\mathcal{N}}$, $\{y_{on}\}_{o \in \mathcal{O}}$, into a integer value variable $\hat{y}_{n} =  \{1,\dots, |\mathcal{O}|\}$ and define $\hat{\textbf{y}} = \{\hat{y}_{n}\}_{n \in \mathcal{N}}$ which is the input to the VAE. 

The VAE architecture consists of two neural networks. The encoder network maps the integer representation $\hat{\textbf{y}}$ onto the parameters of a normal probability distribution $(\mu,\sigma^{2}) = \mathrm{Enc}(\hat{\textbf{y}}; \theta_{e})$ with $\mu \in \mathbb{R}^{N_{h}}$ and $\sigma^{2} \in \mathbb{R}^{N_{h}}$, $N_{h}$ is the dimension of the latent space, and $\theta_{e}$ are the parameters of the encoder. The latent vector is sampled via the reparameterization trick, $z = \mu + \sigma \odot \varepsilon$, with $\varepsilon \sim \mathcal{N}(0,I)$. The decoder performs the inverse mapping, from the continuous to the integer-variable representation of the iterative algorithm, i.e., $\hat{\textbf{y}}^{\mathrm{pred}} = \mathrm{Dec}(z;\theta_{d})$, with $\theta_{d}$ being the learnable parameters of the decoder network. 

To train the VAE, we use a dataset of valid expression trees, i.e., trees whose binary assignments satisfy the constraints presented in Eq.~\ref{eq: one operator}, \ref{eq: operator assignment 1}, \ref{eq: operator assignment 2}. The training seeks the VAE parameters,  $\theta_{e}$ and $\theta_{d}$, that minimize the following loss function
\begin{equation}
    \mathcal{L}_{\text{total}} = \mathcal{L}_{\text{recon}} + \beta \ \mathcal{L}_{\text{KL}}.
\end{equation}
The first term is the reconstruction loss and is equal to the cross-entropy between the original integer representation of the tree $\hat{\textbf{y}}$ and the predicted one $\hat{\textbf{y}}^{\mathrm{pred}}$. This term ensures that the decoder correctly reconstructs the input tree. The second term computes the Kullback-Leibler (KL) divergence, which measures the statistical distance between the encoder's learned distribution $\mathcal{N}(\mu, \sigma^2)$ and a standard normal prior $\mathcal{N}(0, I)$. This term improves the smoothness of the latent space $Z$. The scaling hyperparameter $\beta$ controls the trade-off between these two competing objectives.

\subsection{Evaluating the computational performance of an update function}
For a fixed functional form of the update function, $\hat{\textbf{y}}$, we must compute the values of the constants that optimize the performance function. This corresponds to the following nonlinear optimization problem
\begin{equation} \label{eq: algorithm evaluation}
    \begin{aligned}
m(\hat{\textbf{y}}):= \min_{c} \ & \frac{1}{\mathrm{N_{F}} \mathrm{N_{x_{0}}}} \sum_{f=1}^{\mathrm{N_{F}}} \sum_{p=1}^{\mathrm{N_{x_{0}}}} \phi \big(\hat{\textbf{y}},c,F_{f}, \textbf{x}_{0}^{fp}, \{\textbf{x}_{fpi}\}_{i=1}^{\mathrm{N_{it}}}   \big)\\
\text{s.t.} \ \ & x_{fpki} = x_{fpk,i-1} + g(\hat{\textbf{y}},c,x_{fpk,i-1}, \partial_{x_k} F(\mathbf{x}_{fp,i-1})) \\
& \ \forall k \in \mathcal{N}_x, \forall i \in \mathcal{I}, f \in \mathcal{F}, p \in \mathcal{P}\\
    & x^{\text{lb}} \le x_{fpki} \le x^{\text{ub}} \forall k \in \mathcal{N}_x, \forall i \in \mathcal{I}, f \in \mathcal{F}, p \in \mathcal{P}\\
    & c^{\text{lb}} \le c_{n} \le c^{\text{ub}} \ \forall n \in \mathcal{N} 
\end{aligned}
\end{equation}
Depending on the convergence properties of the iterative algorithm and the bounds on the variables, this problem can be infeasible, i.e., the candidate functional form of the update function is not convergent. In such cases, a large constant value is returned as the computational performance of the candidate algorithm. 

\subsection{Algorithm Discovery via Latent-Space Bayesian Optimization}
The VAE presented above maps the functional form of an update function, parameterized by $\hat{\textbf{y}}$, onto a continuous latent representation $z$, and for a given latent representation $z$, we can obtain the computational performance of the algorithm $\hat{\textbf{y}}=\mathrm{Dec}(z; \theta_{d} )$, $m(\hat{\textbf{y}})=\mathrm{Dec}(z; \theta_{d} )=m(z)$, by solving the problem in Eq.~\ref{eq: algorithm evaluation}.

We use the latent representation to search for new update functions using Bayesian Optimization (BO). BO relies on a Gaussian Process (GP) to approximate the unknown performance function $m(z)$. The BO algorithm is initialized by sampling a random set of valid symbolic trees $\{\hat{\textbf{y}}_{i}\}$ and evaluating the expected computational performance $\{m(\hat{\textbf{y}_{i}})\}$. We use a GP with a Matern kernel as the surrogate model and use the log expected improvement as the acquisition function.

In each BO iteration, the acquisition function is optimized to identify the most promising candidate coordinate $z \in Z$. This coordinate is passed through the decoder, where an \verb|argmax| operation converts the output probabilities into an expression tree, $\hat{\textbf{y}} = \mathrm{Dec}(z;\theta_{d})$. If the decoded tree is invalid, i.e., it does not follow arity rules, we skip the expensive evaluation step and directly assign a large constant penalty to the objective value $m(\hat{\textbf{y}})$. If the tree is valid and the problem in Equation~\ref{eq: algorithm evaluation} is feasible, the objective function value is recorded. This new data point updates the GP, and the loop continues until the predefined number of BO iterations is completed.

\begin{remark}
    \normalfont \textbf{Validity of the reconstructed expression tree:} For a given point in the latent space $z \in Z$, the reconstructed expression tree is not guaranteed to follow arity rules, i.e., be a valid expression tree. This is a common issue in latent-space BO approaches that utilize VAEs, and several approaches have been proposed to address it \cite{kusner2017grammar}. In this paper, we do not use any of these methods since, in our experiments, all the trees suggested by BO followed the arity rules. 
\end{remark}

\section{Results}
To evaluate the proposed framework, we consider two case studies. For each case study, we consider two functions, $\rm{N_{F}}=2$, and two starting points, $\rm{N_{p}}=2$, per function. The maximum number of iterations was set to $\rm{N_{it}}=5$ and the performance function $\phi$ was the norm of the gradient at the last iteration. We compare the computational performance of the discovered algorithms with steepest gradient descent (GD). The step size of GD is tuned by solving an optimization problem that minimizes the number of iterations required to reach a predetermined convergence threshold.  

\subsection{Search space} For both case studies, the update function is modeled as a perfect binary tree of depth $d=2$, i.e., the number of nodes is seven. The binary operator set is $\mathcal{B} = \{+, -, *\}$ and the operand set is $\mathcal{L} = \{x_k, \partial_k f, \text{cst}\}$.

\subsection{VAE architecture and training}
\begin{figure}[t]
    \centering
    \includegraphics[width=\columnwidth]{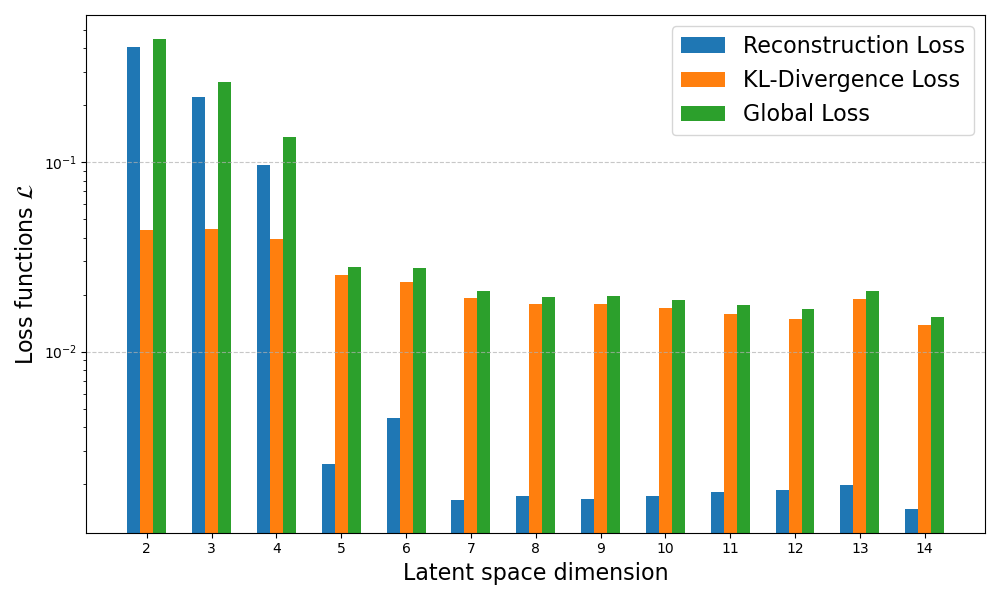}
    \caption{Impact of the latent space dimension on the loss function $\mathcal{L}$ and its components: the reconstruction $\mathcal{L}_{\text{rec}}$ and the KL-divergence $\mathcal{L}_{\text{KL}}$ loss.}
    \label{fig:latent_loss_impact}
\end{figure}

We generate 5000 valid expression trees at random and use 4000 to train the VAE. We use ReLu as an activation function for both the encoder and decoder networks. The encoder processes $\hat{\textbf{y}}$ into 28 neurons, followed by a 64-neuron hidden layer that maps to the latent space. The decoder mirrors this, using a 64- and 28-neuron hidden layer followed by a final output layer. The training was performed using Adam \cite{kingma2014adam} optimizer with a learning rate of $10^{-2}$ for 500 iterations. We performed a sensitivity analysis to identify the latent space dimension by varying the hidden dimension from $2$ to $14$. The loss function and its individual terms are presented in Figure~\ref{fig:latent_loss_impact}. Based on these results, the hidden dimension was set to seven, as increasing it did not significantly improve the loss. Although the VAE training occurs prior to the BO search, the computational overhead is negligible, as training takes approximately $1.66$ seconds. Thus, we excluded this training time from the total reported CPU times below.

\begin{figure*}
    \centering
    \includegraphics[width=\textwidth]{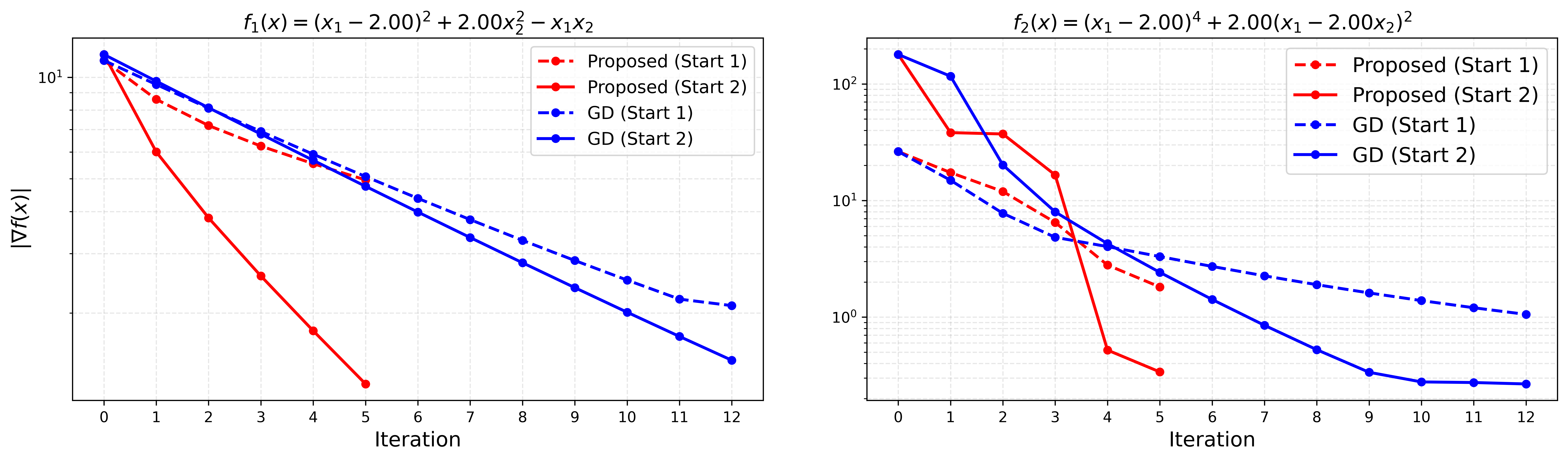}
    \caption{Evolution of gradient norm across iterations for GD and the discovered algorithms in case study 1.}
    \label{fig:nonconvex_trajectories case}
\end{figure*}
\subsection{Update function evaluation}
For the evaluation of the computational performance of a given update function, i.e., solving the problem in Eq.~\ref{eq: algorithm evaluation}, we set the bounds of the variables denoting the constants to $c^{\rm{lb}}=-10, c^{\rm{ub}}=10$ and solve it to global optimality. The model is implemented in Pyomo \cite{hart2017pyomo} and solved using SCIP v.9.2.3 \cite{bolusani2024scip}.

\subsection{Case Study 1}
In the first case study, we consider the minimization of the following functions
\begin{equation}
    \begin{aligned}
        F_1(x) & = (x_1-2.00)^2 + 2.00x_2^2 - x_1 x_2\\
        F_2(x) & = (x_1-2.00)^4 + 2.00(x_1-2.00x_2)^2.
    \end{aligned}
\end{equation}
For each function, we generate two starting points using the NumPy package with a fixed random seed of 1234. The BO algorithm was initialized with 10 random valid trees. After 21 iterations (31 function evaluations in total), the following update rule was discovered
\begin{equation}
    {x}_{i+1} = {x}_i - 0.035 (2.18 + {x}_i) \nabla F({x}_i).
\end{equation}
The total computing time was 121 seconds; 112 seconds were required to evaluate the proposed update function, and 9 seconds to identify the next candidate update function. This update function was not part of the initial set of trees used to initialize BO. 

Unlike standard GD, which utilizes a fixed constant step size, the discovered rule introduces a multiplicative, state-dependent scaling factor $(2.18 + x_i)$. As presented in Fig.~\ref{fig:nonconvex_trajectories case}, this modification yielded significant improvement in the norm of the gradient at the initial iterations. For the first function, the identified algorithm achieved a $L_{1}$ norm for the gradient in five iterations, which can be achieved with steepest gradient descent in 12. Furthermore, for the second function, GD exhibited slow, monotonic convergence. In contrast, the state-dependent multiplier in the discovered rule leads to faster reduction in the $L_{1}$ norm of the gradient. Finally, we note that during the BO iterations, the steepest descent was rediscovered. The ability of the proposed approach to rediscover GD and discover a new variable-step-size update rule shows that the VAE learns a meaningful continuous representation of the space of update functions that BO can explore to discover new algorithms.

\subsubsection{Comparison with mathematical programming approaches}
To evaluate the computational efficiency of the proposed algorithm discovery approach, we compared it with the deterministic MINLP formulation from \cite{mitrai2025automated}. The MINLP model has 2129 constraints, 478 variables (30 binary) and we solve it with BARON v.25.8.5 \cite{zhang2025solving} and SCIP \cite{bolusani2024scip}. After 3600 seconds, the best algorithm found with BARON has an objective value of $6.82$, whereas the solution found with the proposed BO approach after 121 seconds has an objective value of $2.08$.  SCIP identifies a better solution than BARON, with an objective value equal to $2.8$, which corresponds to the steepest gradient descent. In this case study, the results highlight the computational advantage of the proposed approach over existing global optimization solvers.

\subsection{Case Study 2}
In the second case study, we consider the search for new algorithms to solve the linear regression optimization problem, $\min 0.5 \|Ax-b\|^2$. We considered two randomly generated problems with ($A_1$, $b_{1}$) and ($A_2$, $b_{2}$) (using NumPy with random seed of 1234). The BO algorithm was initialized with 20 randomly selected valid trees and then ran for 25 iterations. The following update rule was found 
\begin{equation}
    x_{i+1} = x_i - 0.11 \nabla F   (x_i) - 1.61 \times 10^{-3} x_i.
\end{equation}
The total computational time was $176$ seconds: $157$ seconds were spent evaluating the proposed update functions, and $18$ seconds were spent selecting the next update function. 
\begin{figure}[t]
    \centering
    \includegraphics[width=\columnwidth]{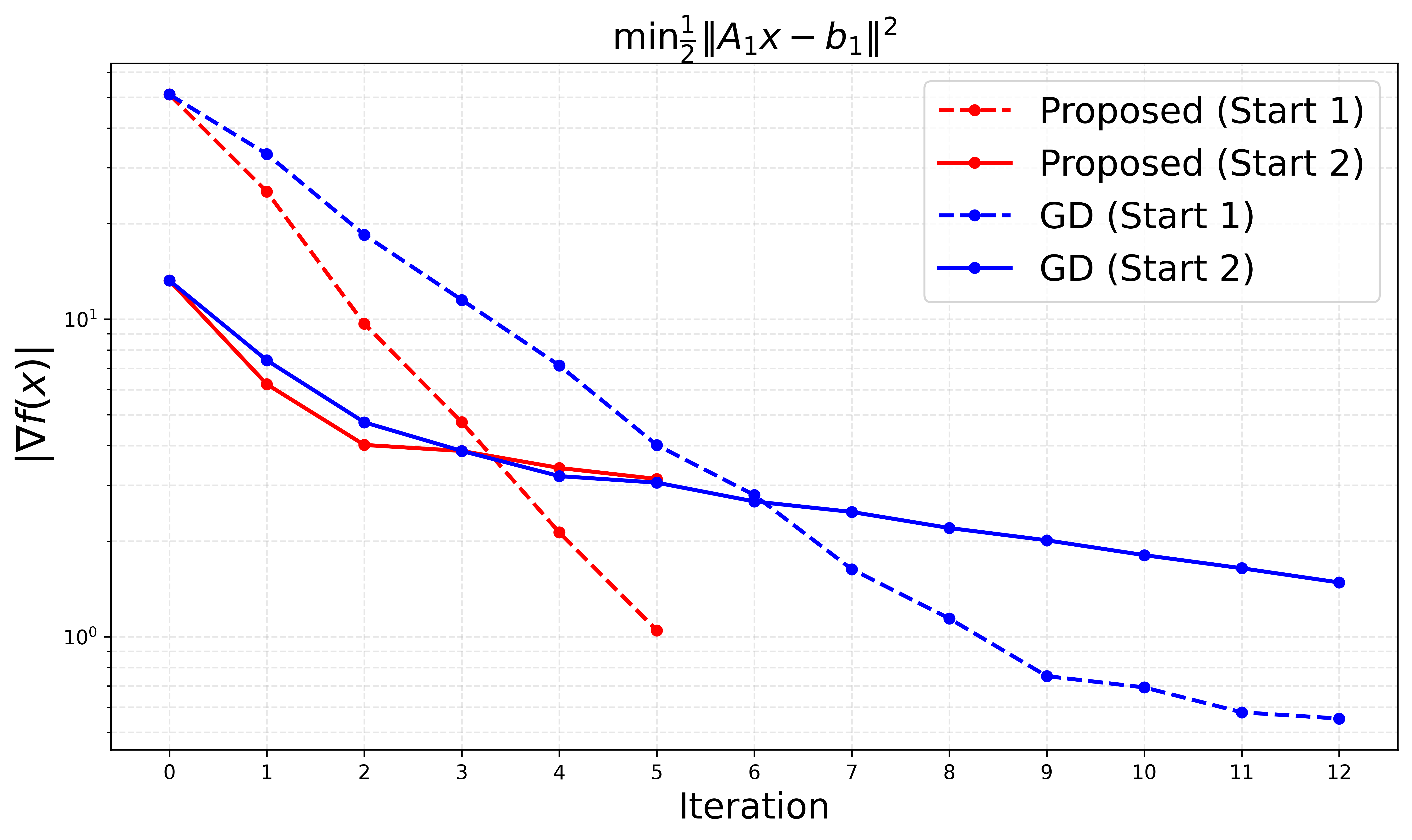}
    \caption{Evolution of gradient norm across iterations for GD and the discovered algorithms for the first system in case study 2.}
    \label{fig:nonconvex_trajectories case 2}
\end{figure}
The two functions differ in the condition number $\kappa$ of the A matrix: for the first, $\kappa_1=16.87$, whereas for the second, $\kappa_2=7.88$. The evolution of the gradient norm for the first function is presented in Fig.~\ref{fig:nonconvex_trajectories case 2}, where we observe that the value of $\|\nabla F(x_{\rm{N_{it}}})\|$ is lower than GD for the first starting point. Overall, the objective of $m(\hat{\textbf{y}})$ with the discovered algorithm is 1.657, whereas with GD 1.679. We note that the search for the update function seeks to minimize the expected computational performance of the iterative algorithm. Thus, for some functions and starting points, the discovered update rule can perform similarly or worse than existing algorithms, such as for the second starting point for the second linear regression problem. 

\section{Conclusion}
In this work, we propose a latent-space Bayesian Optimization approach to discover update rules. We represent the update function as an expression tree, which enables the systematic construction of a superstructure of update functions. To efficiently navigate this discrete space, first, we learn a continuous representation using variational autoencoders, which transforms the algorithm discovery task into a continuous problem. This continuous representation is subsequently used to search for new update functions using BO.

The ability of the proposed approach to discover new update rules was demonstrated through two case studies. In the first case study on non-convex functions, the framework discovered an update rule with a state-dependent step-size. For the second case study, the discovered algorithm led to better computational performance of the function with the higher condition number. In all cases, the computational effort to discover the update rules was significantly lower than that of mathematical programming approaches. The framework is application agnostic and can be used to discover iterative algorithms for general computational tasks. Finally, the computational time can be reduced by solving the evaluation task with local optimization methods.

\section{Acknowledgements}
Financial support from the McKetta Department of Chemical Engineering is gratefully acknowledged.

\bibliographystyle{elsarticle-num}
\bibliography{refs}

\end{document}